\documentclass[12pt]{article}

\usepackage{amsmath,amsthm,amsfonts,amssymb,latexsym,amscd}
\usepackage{pictex,rotating,multirow}

\begin{document}

\newtheorem*{theo}{Theorem}
\newtheorem*{pro}{Proposition}
\newtheorem*{cor}{Corollary}
\newtheorem*{lem}{Lemma}
\newtheorem{theorem}{Theorem}[section]
\newtheorem{corollary}[theorem]{Corollary}
\newtheorem{lemma}[theorem]{Lemma}
\newtheorem{proposition}[theorem]{Proposition}
\newtheorem{conjecture}[theorem]{Conjecture}
\newtheorem{definition}[theorem]{Definition}
\newtheorem{problem}[theorem]{Problem}
\newtheorem{remark}[theorem]{Remark}
\newtheorem{example}[theorem]{Example}
\newcommand{\Naturali}{{\mathbb{N}}}
\newcommand{\Reali}{{\mathbb{R}}}
\newcommand{\Complessi}{{\mathbb{C}}}
\newcommand{\Toro}{{\mathbb{T}}}
\newcommand{\Relativi}{{\mathbb{Z}}}
\newcommand{\HH}{\mathfrak H}
\newcommand{\KK}{\mathfrak K}
\newcommand{\LL}{\mathfrak L}
\newcommand{\as}{\ast_{\sigma}}
\newcommand{\tn}{\vert\hspace{-.3mm}\vert\hspace{-.3mm}\vert}
\def\mA{{\mathfrak A}}
\def\A{{\mathcal A}}
\def\mB{{\mathfrak B}}
\def\B{{\mathcal B}}
\def\C{{\mathcal C}}
\def\D{{\mathcal D}}
\def\F{{\mathcal F}}
\def\H{{\mathcal H}}
\def\J{{\mathcal J}}
\def\K{{\mathcal K}}
\def\L{{\mathcal L}}
\def\N{{\mathcal N}}
\def\M{{\mathcal M}}
\def\O{{\mathcal O}}
\def\P{{\mathcal P}}
\def\S{{\mathcal S}}
\def\T{{\mathcal T}}
\def\U{{\mathcal U}}
\def\W{{\mathcal W}}
\def\Z{{\mathcal Z}}
\def\b{\lambda_B(P}
\def\j{\lambda_J(P}
\def\id{\operatorname{id}}
\def\Ad{\operatorname{Ad}}
\def\aut{\operatorname{Aut}}
\def\SS{{\mathcal S}}
\def\Cb{{\mathbb C}}

\title{Automorphisms of the Cuntz algebras}

\author{Roberto Conti, Wojciech Szyma{\'n}ski}

\date{July 18, 2011}
\maketitle

\renewcommand{\sectionmark}[1]{}

\vspace{7mm}
\begin{abstract}
We survey recent results on endomorphisms and especially on automorphisms of the Cuntz algebras, 
with a special emphasis on the structure of the Weyl group. We discuss endomorphisms globally 
preserving the diagonal MASA and their corresponding actions. In particular, we investigate those 
endomorphisms of $\O_n$ which restrict to automorphisms of the diagonal. We review a combinatorial 
approach to the study of permutative endomorphisms.
All the presented material is put in context with current research topics.
\end{abstract}

\vfill

\vspace{3mm}

\newpage

\section{Introduction}
The Cuntz algebras $\O_n$ form a very well-known class of $C^*$-algebras
which lie at the heart of several interesting problems connecting operator algebras
to quite a number of distinct subjects. In this respect, different properties of $\O_n$ 
turn out to become relevant in different situations and an in-depth analysis of the basic 
general properties of $\O_n$ is absolutely crucial.

In this paper, we present a short overview of the recent progress in the study of the so-called localized 
endomorphisms of the Cuntz algebras $\O_n$,
with special emphasis on the related concepts of full and reduced Weyl groups. Most of the recent results 
described below can be found in \cite{C,CHS1,CHS2,CKS,CRS,CS,HSS,S}. 
Inspired by original Cuntz's insight \cite{Cun2}, we have started a thourough investigation of the
endomorphisms of the Cuntz algebras leaving the canonical MASA $\D_n$ globally invariant.
Automorphisms of $\O_n$ in this family, modulo those fixing $\D_n$ pointwise,
are precisely the elements of the so-called Weyl group of $\O_n$, while those leaving globally invariant 
also the UHF core subalgebra $\F_n$ of $\O_n$ give rise to the reduced Weyl group.
The nontrivial elements of the latter group can be represented by automorphisms that are ``localized'', 
in the sense that they map finite matrix algebras inside $\F_n$ into (larger, in general) finite matrix algebras.
Detailed combinatorial properties of such automorphisms and group-theoretical structural results
for these Weyl groups have been obtained, showing in particular an interesting connection
with the theory of shifts. Several directions for further research have been pointed out, 
in which methods of operator algebras, group theory, symbolic dynamics and combinatorics can hardly 
be separated. Without any pretense to be exhaustive, in this paper we have tried to provide a clue
of the richness of this territory. The study of analogous problems in the more general framework 
of graph algebras is presently under way, with the first set of fundamental results reported in \cite{CHS3}.

\section{Cuntz algebras}
Needless to say, isometries play an important role in functional analysis.
In operator algebras, Coburn investigated the $C^*$-algebra 
$C^*(S)$ generated by a proper isometry $S$, \cite{Cob}.
Some time later, Cuntz introduced and investigated (\cite{Cun1,Cun3}) 
the structure of certain $C^*$-algebras
that carry his name ever since, namely $C^*$-algebras
$\O_n = C^*(S_1,\ldots,S_n)$, $2 \leq n < \infty$, where
$$S_i^* S_i = 1$$
$$\sum_{i=1}^n S_i S_i^* = 1 \ . $$
In particular, $S_i^* S_j = \delta_{ij}1$, for all $1 \leq i,j \leq n$
and
$H := {\rm span} \{ S_1,\ldots,S_n \}$ is a Hilbert space of dimension $n$ in $\O_n$
with support 1.

\medskip
It is useful to introduce some notation and terminology.
Define
$W_n^k := \{1,\ldots,n\}^k$ ($k \geq 1$), the set of $k$-tuples $\alpha = (\alpha_1,\ldots,\alpha_k)$,
where $\alpha_i \in \{1,\ldots,n\}$ for all $1 \leq i \leq k$
and
$W_n := \bigcup_{k=0}^\infty W_n^k$, $W_n^0 = \{0\}$.
Elements of $W_n$ are called multi-indices.
If $\alpha \in W_n^k$, then $l(\alpha) := k$, the length of the word $\alpha$ in the alphabet $\{1,\ldots,n\}$.
Given $\alpha = (\alpha_1,\ldots,\alpha_k) \in W_n$ let
$$S_\alpha := S_{\alpha_1} \ldots S_{\alpha_k}$$ 
($S_0 = 1$ by convention).
The following statement collects some well-known facts about $\O_n$.
\begin{itemize}
\item Every word in $\{S_i, S_i^* \ | \ i = 1,\ldots,n\}$ can be uniquely expressed as
a normal ordered monomial
$S_\alpha S_\beta^*$, for $\alpha, \beta \in W_n$.
\item
$\O_n$ is the closed linear span of $S_\alpha S_\beta^*$,
for $\alpha,\beta \in W_n$.
\item
$\O_n$ is unital, simple, separable, nuclear and purely infinite $C^*$-algebra
with 
$$K_0(\O_n) \cong {\mathbb Z}_{n-1} \quad \text{and} \quad K_1(\O_n) = 0  $$ 
(in particular, the $K$-theory of $\O_2$ is trivial).
\end{itemize}
Furthermore, $\O_n$ is not of type $I$, nor is it even an inductive limit of type $I$ algebras. 
Last but not least, $\O_n$ fits very naturally with several interesting classes of $C^*$-algebras, in particular
it is a Doplicher-Roberts algebra, \cite{DoRo}, a Cuntz-Krieger algebra, \cite{CK}, a graph algebra, \cite{Rae}, 
and a Pimsner algebra, \cite{Pim}.

In operator algebras, Cuntz algebras play an important role
in the general classification program of $C^*$-algebras and group actions, 
structure of crossed products by endomorphisms,
as well as in sector theory and subfactors/index theory, entropy, 
dynamical systems, coding, self-similar sets, 
wavelets (signal processing), 
quantum field theory, abstract group duality,
twisted cyclic cocycles and Fredholm theory.

\subsection{ The Subalgebras of $\O_n$}
Being simple, $\O_n$ has no nontrivial closed two-sided ideals.
However, there are interesting unital $C^*$-subalgebras that will be important in our discussion.
Indeed one has
$$\O_n \supset \F_n \supset \D_n $$
where $\F_n$ is the so-called core UHF-subalgebra 
and $\D_n$ is the diagonal, a canonical MASA with Cantor spectrum.

More in detail, define
$\F_n^k := {\rm span}\{S_\alpha S_\beta^*, l(\alpha)=l(\beta)=k\}$ for $k \geq 1$ 
($\F_n^0 = {\mathbb C}$)
so that $\F_n^k \subset \F_n^{k+1}$ for all $k$.
One has $\F_n^k \simeq M_{n^k} \simeq M_n \otimes \cdots \otimes M_n$ ($k$ factors),
compatible with embeddings $x \mapsto x \otimes 1$.
Set $\F_n := \Big(\bigcup_k \F_n^k \Big)^= \simeq \bigotimes_{k=1}^\infty M_n$,
the UHF-algebra of Glimm type $n^\infty$, with unique trace $\tau$.
In particular, $\F_2 \simeq \bigotimes_{i=1}^\infty M_2$ is the CAR algebra.
There is a faithful conditional expectation $E: \O_n \to \F_n$,
obtained by averaging over the canonical gauge action of $\mathbb T$ defined below,
indeed $\F_n = \O_n^{\mathbb T}$.

$\D_n$ is the commutative $C^*$-subalgebra of $\O_n$
generated by projections $P_\alpha:=S_\alpha S_\alpha^*$, $\alpha \in W_n$;
it is a regular MASA, both in $\F_n$ and $\O_n$.
One has $\D_n \simeq C(X_n)$, where
$X_n := \prod^{\mathbb N} \{1,\ldots,n\}$ 
is the set of all infinite words in the alphabet $\{1,\ldots,n\}$
equipped with the product topology.
Under this isomorphism, $P_\alpha$ corresponds to $\chi_{\sigma_\alpha}$, the 
characteristic function of the cylindrical set of sequences starting with $\alpha$.
The Gelfand spectrum $X_n$ is a Cantor set, 
i.e. a compact, metrizable, totally disconnected space with no isolated points.
Finally, set $\D_n^k := \D_n \cap \F_n^k$, generated by projections 
$P_\alpha$ with $\alpha \in W_n^k$,
isomorphic to the diagonal matrices in $M_{n^k}$.
Then $\D_n$ is the norm closure of $ \Big(\bigcup_k \D_n^k \Big)$.
There exists a faithful conditional expectation from $\F_n$ onto $\D_n$ and whence
from $\O_n$ onto $\D_n$ as well.

\subsection{Endomorphisms}
Let $\U(\O_n)$ be the unitary group of $\O_n$ and
${\rm End}(\O_n)$ be the semigroup of unital $*$-endomorphisms of $\O_n$.
Elements of ${\rm End}(\O_n)$ are automatically injective.
It is a well-known fact that there exists a one-to-one correspondence 
between elements of $\U(\O_n)$ and of ${\rm End}(\O_n)$, 
denoted 
$$U \mapsto \lambda_U \ , $$ 
where $\lambda_U$ is determined by
$$\lambda_U(S_i)=US_i, \quad i=1,\ldots,n \ . $$
The correspondence $U \mapsto \lambda_U$ is not a semigroup morphism, 
rather one has the ``fusion rules'' 
$$\lambda_U \lambda_V = \lambda_{\lambda_U(V)U} \ . $$
For instance, the canonical endomorphism $\varphi: \O_n \to \O_n$,
$$\varphi(a) := \sum_{i=1}^n S_i a S_i^*$$
which restricts to the unilateral shift $x \mapsto 1 \otimes x$ on the UHF subalgebra $\F_n$
satisfies $\varphi = \lambda_\theta$, where
$$\theta:=\sum_{1 \leq i,j \leq n} S_i S_j S_i^* S_j^*$$ 
is the unitary flip operator in $\O_n$, switching the components of $H^2 \simeq H \otimes H$.

Moreover, for all $x \in \F_n^r$ and $m \geq r$, we have
$$\lambda_U(x) = U_m x U_m^*, $$
where
$$U_m := U \varphi(U) \ldots \varphi^{m-1}(U)$$
satisfies the cocycle relation $U_{m+r} = U_m \varphi^m(U_r)$ for all $m,r$.

Also,
${\mathbb T} \ni z \mapsto \lambda_{z 1} =: \alpha_z$ provides the automorphic 
``gauge'' action of $\mathbb T$
(the rescaled periodic modular automorphisms w.r.t. $\omega = \tau \circ E$).
As stated before, $\O_n^{\alpha} = \F_n$ and $\O_n$ is a ${\mathbb Z}$-graded $C^*$-algebra.

The following concept was introduced by Longo, \cite{Lo1,Lo2}, in analogy with QFT.
\begin{definition}
Say that $\lambda_U$ is ``localized'' (or algebraic) if
$U \in \F_n^k$
for some $k$.
\end{definition}
If $U \in \F_n^k$ then
$$\lambda_U (\F_n^r) \subset \F_n^{r+k-1}, \quad r \in {\mathbb N} \ . $$
Moreover, $\lambda_U \lambda_{z 1} = \lambda_{z 1} \lambda_U = \lambda_{zU}$ 
for all $z \in {\mathbb T} = \{z \in {\mathbb C} \ | \ |z| = 1\}$.
It turns out that localized endomorphisms give rise to subfactors with finite index,
which was perhaps the main reason to introduce this notion.

\section{The Weyl groups of the Cuntz algebras}
As stated above, given an endomorphism of $\O_n$ there is an associated index problem.
From this point of view, it would be better to focus on cases for which ${\rm Ind}(\lambda_U)$ 
is finite and larger than one.
On the other extreme, with different motivation, 
in the sequel we will mainly look at cases where $\lambda_U$ is an automorphism of $\O_n$.

Indeed, Cuntz showed long time ago, \cite{Cun2}, that the automorphism group 
of $\O_n$, ${\rm Aut}(\O_n)$, has features similar to semisimple Lie groups,
and proposed a definition of the Weyl group as the normalizer of an infinite-dimensional 
maximal torus in this context,
leaving open the main question of what could be said about the explicit structure of these Weyl groups.

\subsection{General facts about automorphisms of $\O_n$}
Producing nontrivial automorphisms of $\O_n$ is not a simple task. 
We begin by collecting few easy or known facts.

\begin{itemize}

\item The gauge automorphisms of $\O_n$ are induced by scalar unitaries $U = z 1$, 
$z \in {\mathbb C}$, $|z| = 1$, \cite{Cun2}.

\item More generally, one has the so-called Bogolubov automorphisms, \cite{Eva}. 
If $U \in \F_n^1$, then $\lambda_U \in {\rm Aut}(\O_n)$,
$(\lambda_U)^{-1} = \lambda_{U^*}$ and if $U \neq 1$ then $\lambda_U$ is outer.
For such unitaries, $\lambda_U \lambda_V = \lambda_{UV}$;
this is the quasi-free action of $U(n) \simeq \U(\F_n^1)$ on $\O_n$.

A notable example is Archbold's ``flip-flop'' automorphism $\lambda_F$ of $\O_2$, 
\cite{Arch}, where
$$F = S_1 S_2^* + S_2 S_1^* \in \F_2^1$$

\item 
The endomorphism $\lambda_U$ of $\O_n$ is the inner automorphism ${\rm Ad}(V)$, where
$V \in \U(\O_n)$, if and only if $U = V\varphi(V^*)$.

\item 
$\lambda_U \in {\rm Aut}(\O_n)$ if and only if $U^* \in \lambda_U(\O_n)$.
Therefore, in order to check surjectivity of $\lambda_U$, it is enough to know that a single element,
namely $U$ itself, is in the image.
(However, this statement is somewhat self-referential
and thus non terribly useful in practice.)

\item
It may well be the case that $\lambda_U$ is an automorphism but $\lambda_{U^*}$ is not!

\item
One interesting explicit example of an outer automorphisms of $\O_n$ not of Bogolubov-type is provided 
by the Matsumoto-Tomiyama automorphism of $\O_4$, \cite{MaTo}.

\item
All homeomorphisms of the full $n$-shift space $X_n$ 
commuting with the shift transformation extend to automorphisms
of $\O_n$, \cite{CK,Mat1}.

\end{itemize}

In general, there is the problem of finding conditions on a unitary
$U \in {\mathcal U}(\O_n)$ ensuring that $\lambda_U$ is an automorphism
or, better, characterizing all such unitaries and providing methods for 
costructing them.
This is mandatory in order to get more insight into the structure of
${\rm Aut}({\mathcal O}_n)$,
${\rm Out}({\mathcal O}_n)
:= {\rm Aut}({\mathcal O}_n)/{\rm Inn}({\mathcal O}_n)$
and some selected subgroups of these groups.

\bigskip

It is worth to mention a bunch of other abstract results.

\begin{itemize}
\item R\o rdam, \cite{Ror}, and Bratteli-Kishimoto, \cite{BraKi}, consider
\begin{align*}
\U_i & := \{U \in \U(\O_n) \ | \ \lambda_U \in {\rm Inn}(\O_n) \}, \\
\U_a & := \{U \in \U(\O_n) \ | \ \lambda_U \in {\rm Aut}(\O_n) \}, \\
\U_s & := \U(\O_n) \setminus \U_a \ .
\end{align*}
Of course, $\U_i \subset \U_a$. Then
\begin{itemize}
\item $\U_i$ is a dense subset of $\U(\O_n)$,
\item $\U_a$ is a dense $G_\delta$-subset of $\U(\O_n)$,
\item $\U_s$ is a dense $F_\sigma$-subset of $\U(\O_n)$.
\end{itemize}
One proof of this result relies on the so-called Rohlin property for the shift.


\item
${\rm Aut}(\O_n)$ is dense in ${\rm End}(\O_n)$,
in the topology of pointwise norm convergence.

\item It holds
$\overline{\rm Inn}(\O_2) = {\rm Aut}(\O_2)$.
Any two automorphisms of $\O_2$ are approximately inner equivalent,
i.e. there exists a sequence of unitaries $(u_n)_n$ such that $\beta = \lim_n {\rm Ad}(u_n) \circ \alpha$.

\item 
Bratteli-Kishimoto (\cite{BraKi}): 
if $\omega_1$ and $\omega_2$ are two pure gauge-invariant states of the Cuntz 
algebra $\O_n$, 
there is an automorphism $\alpha$ of $\O_n$ such that 
$\omega_1 = \omega_2 \circ \alpha$.

\item
Kishimoto, Ozawa, Sakai (\cite{KiOSa}):
the pure state space is homogeneous under the action of the 
automorphism group (or the subgroup of asymptotically inner automorphisms) for all 
the separable simple $C^*$-algebras
(the asymptotically inner automorphisms are a normal subgroup of approximately inner 
automorphisms).
\end{itemize}

\subsection{The Weyl groups}

We systematically employ the following notation.
If $\mathfrak A$ is a unital $C^*$-algebra, denote by
${\rm Aut}({\mathfrak A})$ and 
${\rm Inn}({\mathfrak A})$ the group of $*$-automorphisms
and its normal subgroup of inner automorphisms, respectively,
and we denote by ${\rm Out}({\mathfrak A}) := {\rm Aut}({\mathfrak A})/{\rm Inn}({\mathfrak A})$
the corresponding quotient with canonical projection
$\pi: {\rm Aut}({\mathfrak A}) \to {\rm Out}(\mathfrak A)$.

If
${\mathfrak B} \subset {\mathfrak A}$ is a unital $C^*$-subalgebra, then we define
${\rm Aut}({\mathfrak A},{\mathfrak B}) = \{\alpha \in {\rm Aut}({\mathfrak A}) \ | \ \alpha({\mathfrak B}) = {\mathfrak B}\}$
and
${\rm Aut}_{\mathfrak B}({\mathfrak A}) =
\{\alpha \in {\rm Aut}({\mathfrak A}) \ | \ \alpha(b)=b, \forall b \in {\mathfrak B}\}$.
Clearly
${\rm Aut}_{\mathfrak B}({\mathfrak A}) \triangleleft {\rm Aut}({\mathfrak A},{\mathfrak B})
\subseteq N_{{\rm Aut}({\mathfrak A})}({\rm Aut}_{\mathfrak B}({\mathfrak A}))$.
Denoting by
$N_{\mathfrak A}({\mathfrak B}) = \{u \in \U({\mathfrak A}) \ | \ u {\mathfrak B} u^* = {\mathfrak B} \}$ 
the (unitary) normalizer of $\mathfrak B$ in $\mathfrak A$, there are natural maps
$N_{\mathfrak A}({\mathfrak B}) \to {\rm Aut}({\mathfrak A},{\mathfrak B}) \cap {\rm Inn}({\mathfrak A})$
and ${\rm Aut}({\mathfrak A},{\mathfrak B}) \to {\rm Aut}({\mathfrak B})$.

\medskip

The following result by Cuntz is at the basis of the definition of the (reduced and full) Weyl groups.
Notice that in our context the (unitary) normalizers satisfy
$N_{\F_n^k}(\D_n^k) 
\subset N_{\F_n^{k+1}}(\D_n^{k+1}) \subset \ldots 
\subset N_{\F_n}(\D_n)$.

\begin{theorem} (Cuntz, \cite{Cun2})
\begin{itemize}
\item
$\U(\D_n) \simeq \lambda(\U(\D_n)) = \lambda(\U(\D_n))^{-1} $ 

\item
$\O_n^{\lambda(\U(\D_n))} = \D_n$,
${\rm Aut}_{\D_n}(\O_n) = \lambda(\U(\D_n))$ 

\item
$N_{{\rm Aut}(\O_n)}(\lambda(\U(\D_n)))
= {\rm Aut}(\O_n,\D_n)
= \lambda(N_{\O_n}(\D_n))^{-1}$

\item
${\rm Aut}(\O_n,\D_n) \cap {\rm Aut}(\O_n,\F_n)
= \lambda(N_{\F_n}(\D_n))^{-1}$
\end{itemize}
\end{theorem}

Here, for $E \subset \U(\O_n)$, we denote 
$$\lambda(E)^{-1} := \{\lambda_U \ | \ U \in E \} \cap {\rm Aut}(\O_n) \ . $$

Thus, $\lambda(\U(\D_n))$ plays the role of a ``maximal torus'', 
i.e. a maximal abelian subgroup which is limit of finite-dimensional tori,
the second line states that the pair $\D_n \subset \O_n$ is Galois-closed,
while the third and fourth lines provide, after taking the quotient with ${\rm Aut}_{\D_n}(\O_n)$, 
the full and reduced Weyl groups, respectively.
\bigskip

The next problem is thus to find out
which elements of $N_{\O_n}(\D_n)$,
resp. of $N_{\F_n}(\D_n)$,
induce automorphisms of $\O_n$.
To this aim, the structure of normalizers is provided by the following result.

\begin{theorem}
(Power, \cite{Pow})
\begin{itemize}

\item
$N_{\O_n}(\D_n) = {\mathcal S}_n \cdot \U(\D_n)$, where
$\S_n : =
\{u \in \U(\O_n) \ | \ u = \sum'_k S_{\alpha_k} S_{\beta_k}^* \}$ is
the subgroup of $\U(\O_n)$ of unitaries that can be written as finite sum of words in $S_i$ and $S_i^*$;

\item
$N_{\F_n}(\D_n) = \Big(\bigcup_k N_{\F_n^k}(\D_n^k)\Big) \cdot \U(\D_n)
= \P_n \cdot \U(\D_n)$, where
${\mathcal P}_n =
\{u \in \U(\F_n) \ | \ u = \sum'_k S_{\alpha_k} S_{\beta_k}^*, l(\alpha_k)=l(\beta_k) \ \forall k \}$

\end{itemize}
\end{theorem}

By a result of Nekrashevych, \cite{Nek}, $\S_n$ is isomorphic to the celebrated 
Higman-Thompson group $G_{n,1}$.
In particular, a copy of the Thompson group sits naturally in $\S_2 \subset \U(\O_2)$.
Also, $\P_n$ is an inductive limit of permutation groups ${\mathbb P}_{n^k}$ w.r.t. 
strictly diagonal embeddings.

\bigskip

All in all, 
Cuntz problem boils down to recognizing which unitaries in $\S_n$, resp. $\P_n = \S_n \cap \F_n$, 
induce automorphisms of $\O_n$. We will mainly focus on $\P_n $.

\bigskip

Let $P_n^k$ be the group of permutations of $W_n^k$, clearly isomorphic to ${\mathbb P}_{n^k}$.

To any $\sigma \in P_n^k$ one associates a unitary in $\F_n^k$ via
$$u_\sigma = \sum_{\alpha \in W_n^k} S_{\sigma(\alpha)}S_\alpha^* \ . $$
Then $\sigma \mapsto u_\sigma$ is an isomorphism of $P_n^k$ 
with its image $\P_n^k = \S_n \cap \F_n^k$,
that can be further identified with the set of permutation matrices in $M_{n^k}$,
and $\P_n = \bigcup_k \P_n^k$.

\begin{proposition}
Let $w$ be a unitary in $\O_n$.
\begin{itemize}
\item[(a)]
If $w \in \U(\O_n)$ then 
$\lambda_w(\F_n)=\F_n$ if and only if
$\lambda_w \in {\rm Aut}(\O_n)$ and $w \in \U(\F_n)$;
\item[(b)]
If $\lambda_w \in {\rm Aut}(\O_n)$ then $\lambda_w(\D_n) = \D_n$ if and only
if $w \in N_{\O_n}(\D_n)$.
\item[(c)]
If $\lambda_w(\D_n) = \D_n$ then $\lambda_w$ is an irreducible endomorphism of $\O_n$,
i.e. $\lambda_w(\O_n)' \cap \O_n = {\mathbb C}$.
\end{itemize}
\end{proposition}
As long as we consider endomorphisms of $\O_n$ induced by unitaries $w$ in
$\cup_k \P_n^k \subset N_{\F_n}(\D_n) = N_{\O_n}(\D_n) \cap \F_n$,
when they are automorphisms they also provide, by restriction,
automorphisms of $\D_n$ and $\F_n$;
when they only satisfy the weaker condition $\lambda_w(\D_n)=\D_n$ they still
act irreducibly on $\O_n$.

\begin{theorem} One has

\begin{itemize} 

\item $N_{\O_n}(\D_n) \simeq \U(\D_n) \rtimes \S_n$ (w.r.t. action by conjugation);

\item $N_{\F_n}(\D_n) \simeq \U(\D_n) \rtimes \P_n$;

\item ${\rm Aut}(\O_n,\D_n) \simeq \lambda(\U(\D_n)) \rtimes \lambda(\S_n)^{-1}$.

In particular, $\lambda(\S_n)^{-1}$ is a subgroup of ${\rm Aut}(\O_n,\D_n)$;

\item ${\rm Aut}(\O_n,\D_n) \cap {\rm Aut}(\O_n,\F_n)
\simeq \lambda(\U(\D_n)) \rtimes \lambda(\P_n)^{-1}$;

In particular, $\lambda(\P_n)^{-1}$ is a subgroup of ${\rm Aut}(\O_n,\D_n) \cap {\rm Aut}(\O_n,\F_n)$.

\end{itemize}

\end{theorem}

\subsection{The reduced Weyl group and shift automorphisms}

\begin{proposition}
\label{localizedinverse}
Let $w \in \P_n^k$ and suppose that $\lambda_w \in {\rm Aut}(\O_n)$,
then the inverse $\lambda_w^{-1}$ is also localized. More precisely,
$\lambda_w^{-1}$ is induced by a unitary in  
$\P_n^h$, with $h \leq n^{2(k-1)}$.
\end{proposition}

\begin{proposition}
Let $w \in \P_n$. If $\lambda_w \in {\rm Inn}(\O_n)$ then there exists $z \in \P_n$ such that $w = z \varphi(z^*)$.
Moreover, for $k \geq 2$, if $w \in \P_n^k$ then $z \in \P_n^{k-1}$.
\end{proposition}

\bigskip

There is an isomorphism $\P_n \to \lambda(\P_n)^{-1} \cap {\rm Inn}(\O_n)$, via $u \mapsto {\rm Ad}(u)$.
Thus there is a short exact sequence
$$1 \to \P_n \stackrel {{\rm Ad}} \to \lambda(\P_n)^{-1} \to \pi(\lambda(\P_n)^{-1}) \to 1$$

If $\lambda_u\in\lambda(\P_n)^{-1}$ then $\lambda_u(\P_n)=\P_n$ and thus $\lambda_u$ 
is an automorphism of the locally finite group $\P_n$. It is known that groups $\P_n$ and $\P_m$ 
are isomorphic if and only if $n$ and $m$ have the same prime factors, \cite{KrSu}. 
We note that automorphisms of inductive 
limits of symmetric groups have been studied by group theorists before, e.g. see \cite{LaSh}. 
In fact, we have the following exact sequence 
$$1 \to \P_n \stackrel {{\rm Ad}} \to \lambda(\P_n)^{-1}|_{\P_n} \to \pi(\lambda(\P_n)^{-1}) \to 1$$
Indeed, suppose that $u\in\P_n$ and $\lambda_u(w)=w$ for all $w\in\P_n$. If $w\in\P_n^1$ then this 
means that $uw=wu$. Then $\lambda_u(\varphi(w))=\varphi(w)$ implies $u\varphi(w)=\varphi(w)u$ 
and, by induction, $u\varphi^k(w)=\varphi^k(w)u$ for all $w\in\P_n^1$ and $k=0,1,2,\ldots$ Denote 
$C=(\P_n^1)'\cap\F_n^1$, a 2-dimensional space spanned by $1$ and $\sum_{i,j=1}^n S_iS_j^*$. 
If $u\in\P_n^k$, the above implies that viewing $\F_n^k$ as $M_n\otimes\ldots \otimes M_n$ (a $k$-fold 
tensor product) we have $u\in C\otimes\ldots\otimes C$. If $n\geq3$ then $\P_n^k\cap (C\otimes
\ldots\otimes C)=1$, and thus $u=1$. If $n=2$ this imples that $u=\sigma_1\otimes\ldots\otimes\sigma_k$, 
with $\sigma_j\in\P_2^1$. That is, $u=\sigma_1\varphi(\sigma_2)\ldots\varphi^{k-1}(\sigma_{k-1})$. 
Suppose $\sigma_1\neq\id$, otherwise the argument is analogous. Depending on whether $\sigma_2=\id$ or 
$\sigma_2=\sigma_1$, $\lambda_u$ restricted to $\P_2^2$ acts as $\Ad(\sigma_1\varphi(\sigma_1))$ or 
$\Ad(\sigma_1)$, respectively. In either case there exists $w\in\P_2^2$ such that $\lambda_u(w)\neq w$, 
since $\P_2^2$ being the symmetric group on 4 letters has trivial center.

\medskip
Let $G_n := \lambda(\P_n)^{-1} \subset {\rm Aut}(\O_n)$. 
The following result shortly summarizes some of our achievements, thus sheding some light on the structure of these Weyl groups, and also establishes an intriguing connection with the theory of shifts.

\begin{theorem}
With the above notation, it holds:
\begin{itemize} 

\item $\pi(G_n)$, $n \geq 3$ is non-amenable; actually it contains the modular group 
${\mathbb Z}_2 * {\mathbb Z}_3 (\simeq {\rm SL}_2({\mathbb Z}))$;

\item
$\pi(G_2) $ is non-amenable;
(the class of) Archbold's flip-flop is the simplest nontrivial element,
the next one being induced by a permutative unitary in $\P_2^4$;

\item $\pi(G_n)$ is residually finite, for every $n \geq 2$.

\end{itemize}
\end{theorem}

Loosely speaking, the general method for proving these 
results is summarized by the following diagram:
$$
\begin{array}{ccccccccc}
& & {\rm Aut}(X_n) & & & & & & \\
& & \cap & & & & & & \\
G_n
& \stackrel{\simeq}\to & {\mathfrak G}_n & \to & {\mathfrak G}_n/I{\mathfrak G_n} \simeq \pi(G_n)
& \hookrightarrow & {\rm Aut}(\Sigma_n)/\langle \sigma \rangle \\
\cap & & \cap & & \cap & & & & \\
 {\rm Aut}(\O_n) & & {\rm Aut}(\D_n) & & {\rm Out}(\O_n) & & & & 
\end{array}
$$
where we have 
$$ G_n\cong \frac{{\rm Aut}(\O_n,\D_n) \cap {\rm Aut}(\O_n,\F_n)}{{\rm Aut_{\D_n}(\O_n)}}. $$

Indeed, it turns out that $G_n = \lambda(\P_n)^{-1}$ is isomorphic to the restricted Weyl group of $\O_n$.
Moreover, $G_n$ is isomorphic (via restriction) with the group ${\mathfrak G}_n$ 
of homeomorphisms of the full $n$-shift space $X_n$ which eventually commute, along with 
their inverses, with the shift transformation 
(i.e., they commute with a sufficiently large power of the shift).
For prime $n$, the restricted outer Weyl group $\pi(G_n)$ is isomorhic with the automorphism group of the full two-sided $n$-shift $\Sigma_n$ divided by its center (generated by the shift).

\subsection{The full Weyl group} 

The study of $\lambda(\S_n)^{-1}$ and its quotient modulo inners, i.e. the full outer Weyl group,
are presently under investigation.
We only mention here that one first objective is to find necessary and sufficient conditions for $w \in \S_n$ such that

\begin{itemize}

\item $\lambda_w \in {\rm Aut}(\O_n)$

\item $\lambda_w(\D_n)=\D_n$

\end{itemize}
and moreover to characterize intrinsically the group of homeomorphisms of $X_n$ arising in this way.
In other words, one would like to recognize which automorphisms 
of $\D_n$ can be obtained by restricting automorphisms (or even 
endomorphisms) of $\O_n$. At this point, we can present the following approach to deciding  
which endomorphisms $\lambda_u$, $u\in{\mathcal S}_n$, give rise to automorphisms of $\D_n$. 

\bigskip
For a $u\in\SS_n$ and a $j\in\{1,\ldots,n\}$, let $a_j^u:\D_n\to\D_n$ be a linear map 
defined as 
\begin{equation}\label{aju}
a_j^u(x) = S_j^*u^*xuS_j. 
\end{equation}
(Similar maps were used earlier in the analysis of permutative endomorphisms, \cite{CP,S,CS,HSS}.) 
If no confusion may arise then we write $a_j=a_j^u$. Since $u^*\D_nu=\D_n$, for any $x,y$ we have 
$a_j(xy)=\sum_{i=1}^n S_j^*u^*xuS_iS_i^*u^*yuS_j=a_j(x)a_j(y)$, and the map $a_j$ is in fact 
a unital $*$-homomorphism. We also note that for all $k=0,1,2,\ldots$ we have 
\begin{equation}\label{adtower}
a_j(\lambda_u(\D_n^{k+1}))\subseteq\lambda_u(\D_n^k). 
\end{equation}
Indeed, for any $x_i\in\D_n^{k}$, $i=1,\ldots,n$, we have 
$$ a_j\lambda_u\left(\sum_{i=1}^n P_i\varphi(x_i)\right)=\sum_{i=1}^n(a_j\lambda_u)(P_i)
   (a_j\lambda_u\varphi)(x_i)=a_j\lambda_u\varphi(x_j)=\lambda_u(x_j), $$
and the claim follows. For $\mu\in W_n^k$ we denote 
\begin{equation}\label{Tmu}
T_\mu^u=a_{\mu_k}^u \ldots a_{\mu_2}^u a_{\mu_1}^u. 
\end{equation}
Again, we will usually write $T_\mu=T_\mu^u$, if no danger of confusion arises. One can easily 
verify by induction that for all $k$ and all $d\in\D_n$ we have 
\begin{equation}\label{ukt}
u_k^* d u_k = \sum_{\mu\in W_n^k}S_\mu T_\mu(d) S_\mu^*. 
\end{equation}

\begin{theorem}\label{pm1criterion}
Let $u\in\SS_n$. Then $\lambda_u|_{\D_n}\in\aut(\D_n)$ if and only if there exists an integer 
$m$ such that $T_\mu(\D_n^{\ell'})\subseteq \Cb 1$ for all $\mu\in W_n^m$. 
\end{theorem}
{\em Proof.}
At first suppose that $\lambda_u|_{\D_n}\in\aut(\D_n)$. Let $p$ be a minimal projection in $\D_n^{\ell'}$, 
and let $k$ be so large that $p\in\lambda_u(\D_n^k)$. Then for each $\mu\in W_n^k$ we see that  
$T_\mu(p)$ is a scalar by virtue of formula (\ref{adtower}). Thus it suffices to take $m$ equal to the 
maximum of all such $k$'s. 

Conversely, let $m$ be an integer such that $T_\mu(\D_n^{\ell'})\subseteq 
\Cb 1$ for all $\mu\in W_n^m$. Fix an element $d$ of $\D_n^{\ell'}$. Then $u_m^* d u_m$ belongs 
to $\D_n^m$ by virtue of formula (\ref{ukt}). Thus for each $k>m$ we have 
$$ u_k^* d u_k = \varphi^{k-1}(u^*)\ldots\varphi^m(u^*) (u_m^* d u_m) \varphi^m(u^*)\ldots 
   \varphi^{k-1}(u) = u_m^* d u_m, $$ 
and we see that the sequence $\{u_k^* d u_k\}$ stabilizes from $k=m$. This easily implies that 
element $d$ belongs to $\lambda_u(\D_n)$. Since $d\in\D_n$ was arbitrary,  $\lambda_u|_{\D_n}$ 
is surjective and hence $\lambda_u|_{\D_n}\in\aut(\D_n)$. 
\hfill$\Box$

\medskip
The above theorem does not seem algorithmic, since it is not immediately clear how big integer 
$m$ one needs to consider. However, in certain special cases this problem disappears. Indeed, 
let $u\in\SS_n$ be such that $u=\sum_{(\alpha,\beta)\in\J}S_\alpha S_\beta^*$. We denote 
\begin{equation}\label{ku}
K_u=\max\{||\alpha|-|\beta|| : (\alpha,\beta)\in\J\}. 
\end{equation}
In particular, $u\in\P_n$ if and only if $K_u=0$. Furthermore, if $u\in\SS_n$ is 
such that $K_u\leq 1$ and $\ell'=\ell'(\J)$ then for each $j\in\{1,\ldots,n\}$ we have $a_j^u(\D_n^{\ell'})\subseteq(\D_n^{\ell'})$. This implies that the criterion of Theorem 
\ref{pm1criterion} becomes algorithmic for such unitaries. In fact, this criterion 
can then be checked via a combinatorial procedure similar to the one developed for permutative 
endomorphisms in \cite{CS}. We also note that endomorphisms $\lambda_u$ with $K_u=1$ 
have been considered and played a role earlier in \cite{CRS}. 


\subsection{A few remarks on automorphisms of $\D_n$ obtained from $\O_n$}

During our investigations of the Weyl group of $\O_n$, the following facts have been shown:

\begin{itemize}
\item 
there are proper endomorphisms of $\O_n$ that restrict to automorphisms of $\D_n$.
However, any unital endomorphism of $\O_n$ which fixes the diagonal $\D_n$ point-wise is automatically
surjective, i.e. it is an element of ${\rm Aut}_{\D_n}(\O_n)$;

\item 
the restriction map ${\rm Aut}(\O_n,\D_n) \to {\rm Aut}(\D_n)$ is not surjective 
(and its image is not normal, as ${\rm Aut}(\D_n)\simeq {\rm Homeo}(\C_n)$ is a simple group);

\item
there are product-type automorphisms of $\D_n$
that do not extend to (possibly proper) endomorphisms of $\O_n$;
in case of $\D_2$, consider e.g.
$\bigotimes_{i=1}^\infty {\rm Ad}(u_i)$, where 
$$u_i = \begin{cases} 1 & i \quad \mbox{even} \\ 
( \begin{smallmatrix} 0 & 1 \\ 1 & 0
\end{smallmatrix} ) & i \quad \mbox{odd} \end{cases}$$
and we have realized $\D_2$ as an infinite tensor product over $\mathbb N$ of diagonal matrices of size $2$;

\end{itemize}

\medskip
As another aspect of the subtle interplay between $\O_n$ and $\D_n$ we would like to mention the following
decomposition problem (Cuntz): if $\lambda_U \in {\rm Aut}(\O_n)$, is there a unitary $V$ in $\O_n$ such that
$\lambda_U(\D_n) = {\rm Ad}(V)(\D_n)$?

\bigskip
Even the inner part of the group $\lambda(\SS_n)^{-1}$, that is the group $\SS_n$ itself acting on 
$\O_n$ by $\Ad(u)$, has an extraordinarily complicated structure and action on the diagonal $\D_n$. 
In particular, a natural question arises if there exists a subgroup $\Gamma_n$ of $\SS_n$ such that 
$\O_n = \D_n \rtimes \Gamma_n$? Spielberg gave an affirmative answer to this question in the 
case of an {\em even} $n$, \cite{Spi}. Namely, in the case $n=2$, one takes 
$\Gamma_2 = {\mathbb Z}_2 * {\mathbb Z}_3$, with generators $u$ (order 2) and $t$ (order 3). 
With $\partial \Gamma_2$ the  boundary of $\Gamma$, there exists a natural action of $\Gamma_2$ on 
$\partial \Gamma_2$ such that
$$C(\partial \Gamma_2) \rtimes \Gamma_2 \simeq M_2(\O_2).$$
But $C(\partial \Gamma_2) \simeq \D_2$ and $M_2(\O_2) \simeq \O_2$! In this picture, one may identify 
$u = S_1 S_2^* + S_2 S_1^*$ and $t = S_2 S_1 S_1^* S_1^* + S_2 S_2^* S_1^* + S_1 S_1 S_2^*$. 
The general even case is handled similarly, with $\Gamma_{2m} = {\mathbb Z}_m 
\times ({\mathbb Z}_2 * {\mathbb Z}_{2m+1})$, and one has 
$$ \O_{2m}=\D_{2m}\rtimes\Gamma_{2m}. $$
However, the existence of such a crossed product decomposition in the case $n$ odd is unknown to us. 

\section{Constructive results: permutative automorphisms and trees}

Combining above results with brute force computation 
(with a lot of help from computer) it has been possible to determine all
permutative automorphisms of $\O_n$ at level $k$, i.e. associated to unitaries
in $\P_n^k$, for small values of $n \geq 2$ and $k$.
More precisely, by explicit computations the numbers

\bigskip

$d_n^k := \# \{\sigma \in P_n^k \ | \lambda_\sigma \in {\rm Aut}(\O_n)\} \equiv \# \lambda(\P_n^k)^{-1}$

\bigskip

$b_n^k = \# \{\sigma \in P_n^k \ | \ \lambda_\sigma|_{\D_n} \in {\rm Aut}(\D_n)\} $,

\smallskip
\noindent have been completely determined for $n+k \leq 6$
(besides, notice that $n!^{n^{k-1}} \ \vert \ b_n^k$).

Without digging into details, the complicated combinatorial structure has been dealt with
by making a clever use of ($n$-tuples of) labeled rooted trees, 
with vertices indexed by elements in $W_n^{k-1}$.
Equivalently, one may consider a directed graph in which the set of edges is decomposed into 
$n$ spanning rooted trees (the root being a fixed point). 

In the case of $\P_2^4$, it has been possible also to compute explicitly part of the
structure of the infinite subgroup of ${\rm Out}(\O_2)$ generated by the 14 inner equivalence classes 
which were found.

For $\O_n$ at level $k$ one has to consider $n^k!$ permutations.
A case-by-case brute force computation is unfeasible as it involves some manipulation of large matrices
(with size that could grow up to the order of $n^{n^{2k}}$ or so). Moreover, 
there are simply too many cases to consider as the following figures illustrate:

\medskip

$\O_2$: 

$2! = 2$
 
$2^2! =  24$ 

$2^3! = 40320$ 

$2^4! = 20,922,789,888,000$, ...

\medskip

$\O_3$: 

$3! = 6$ 

$3^2 ! = 362880$ 

$3^3! = 10,888,869,450,418,352,160,768,000,000$, 

$\ldots$

\medskip

Some simplifications are possible, exploiting the action of inner and Bogolubov automorphisms,
but they do not affect significantly the scale of the problem. Need help!
Surprisingly enough, labeled rooted trees $T$ come to the rescue, where
$$\# V(T) = \# E(T) = n^{k-1}$$
Why trees ?
Given $\sigma \in P_n^k$ define functions 
$$f^\sigma_i: W_n^{k-1} \to W_n^{k-1}$$ 
for $i=1,\ldots n$ by
$$f^\sigma_i(\alpha) = \beta : \Leftrightarrow \mbox {there exists $m \in \{1,\ldots,n\}$ such that} \  (i,\alpha)=\sigma(\beta,m) \ . $$
Then a necessary condition for $\sigma$ to give rise to an automorphism is that the ``diagrams'' of all the $f^\sigma_i$'s  are rooted trees, where the root is the unique fixed-point.
In particular, the vertices of these trees are labeled by $W_n^{k-1}$.
Moreover, in order to get automorphisms, 
this labeling must also induce a certain partial order relation on $W_n^{k-1} \times W_n^{k-1}$!

\begin{example}
{\rm The pair of labeled trees corresponding to
$\sigma = {\rm id}$ in $P_2^3$. All the edges are downward oriented. }


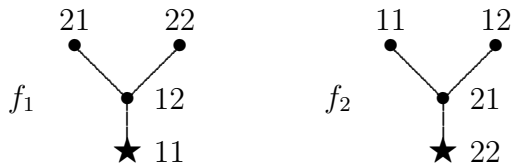
\[ \beginpicture
\setcoordinatesystem units <0.7cm,0.7cm>
\setplotarea x from 4 to 9, y from -1 to 1

\put {$\bullet$} at -1 1
\put {$\bullet$} at 1 1
\put {$\bullet$} at 0 0
\put {$\bigstar$} at 0 -1

\setlinear
\plot -1 1 0 0 /
\plot 1 1 0 0 /
\plot 0 0 0 -1 /
\put {$f_1$} at -2 0
\put {$21$} at -1 1.5
\put {$22$} at 1 1.5
\put {$12$} at 0.8 0
\put {$11$} at 0.8 -1

\put {$\bullet$} at 5 1
\put {$\bullet$} at 7 1
\put {$\bullet$} at 6 0
\put {$\bigstar$} at 6 -1

\setlinear
\plot 5 1 6 0 /
\plot 7 1 6 0 /
\plot 6 0 6 -1 /
\put {$f_2$} at 4 0
\put {$11$} at 5 1.5
\put {$12$} at 7 1.5
\put {$21$} at 6.8 0
\put {$22$} at 6.8 -1

\endpicture \]
\end{example}

\begin{example}
{\rm Let $u\in{\mathcal P}_n^1$, so that $\lambda_u$ is a Bogolubov
automorphism of $\O_n$. If we view $u$ as an element of ${\mathcal P}_n^k$
then all $n$ unlabeled trees corresponding to $u$ are identical; the root
receives $n-1$ edges from other vertices, each other vertex receives either
none or $n$ edges, and the height of the tree (the length of the longest path
ending at the root) is minimal and equal to $k-1$. In particular, all
such unitaries have the corresponding $n$-tuples of unlabeled trees
identical with those of the identity.}
\end{example}

\begin{example}
{\rm 
In the case of $\P_2^4$, we have 
\begin{align*}
& \# \{\lambda_w \ | \ w\in\P_2^4 \mbox{ and } \lambda_w|_{\D_2}
\in{\rm Aut}(\D_2)\} = 8! \cdot 2^8 \cdot 17 = 175,472,640 \ , \\
& \# \{\lambda_w \ | \ w\in\P_2^4 \mbox{ and } \lambda_w
\in{\rm Aut}(\O_2)\} = 8! \cdot 14 = 564,480 \ .
\end{align*}
Thus in $\lambda(\P_2^4)^{-1}$ there are exactly $14$ representatives
of distinct inner equivalence classes.
Moreover, it is possible to verify by explicit calculations that the infinite 
subgroup of ${\rm Out}(\O_2)$ so obtained is generated by involutions.

The following picture shows all rooted trees which may arise from permutative endomorphisms 
$\lambda_u$, $u\in\P_2^4$, 
grouped according to the number $f$ of leaves. They are 23. 
Remarkably, only the two in the bottom row may arise when 
$\lambda_u$ is an automorphism of $\O_2$!

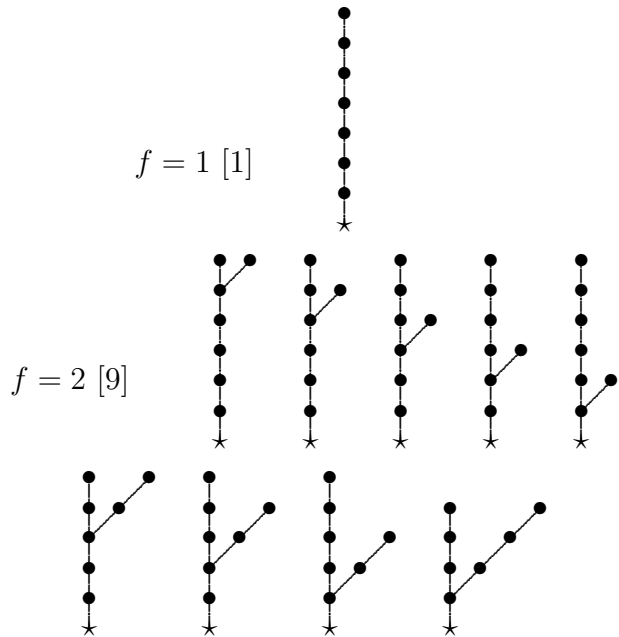
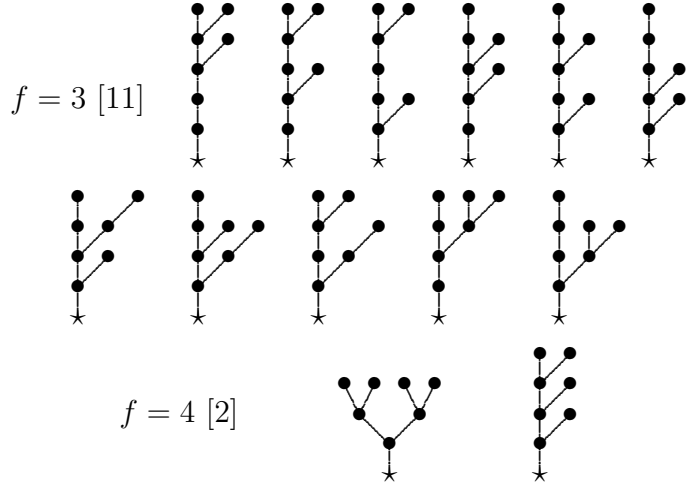
\[ \beginpicture
\setcoordinatesystem units <0.4cm,0.4cm>
\setplotarea x from 0 to 5, y from -1 to 1

\put {$f=1$ [1]} at -5 0

\put {$\bullet$} at 0 5
\put {$\bullet$} at 0 4
\put {$\bullet$} at 0 3
\put {$\bullet$} at 0 2
\put {$\bullet$} at 0 1
\put {$\bullet$} at 0 0
\put {$\bullet$} at 0 -1
\put {$\star$} at 0 -2

\setlinear
\plot 0 5 0 4 /
\plot 0 4 0 3 /
\plot 0 3 0 2 /
\plot 0 2 0 1 /
\plot 0 1 0 0 /
\plot 0 0 0 -1 /
\plot 0 -1 0 -2 /

\endpicture \]

\[ \beginpicture
\setcoordinatesystem units <0.4cm,0.4cm>
\setplotarea x from 0 to 5, y from -1 to 1

\put {$f=2$ [9]} at -5 0

\put {$\bullet$} at 1 4
\put {$\bullet$} at 0 4
\put {$\bullet$} at 0 3
\put {$\bullet$} at 0 2
\put {$\bullet$} at 0 1
\put {$\bullet$} at 0 0
\put {$\bullet$} at 0 -1
\put {$\star$} at 0 -2

\setlinear
\plot 1 4 0 3 /
\plot 0 4 0 3 /
\plot 0 3 0 2 /
\plot 0 2 0 1 /
\plot 0 1 0 0 /
\plot 0 0 0 -1 /
\plot 0 -1 0 -2 /

\put {$\bullet$} at 4 3
\put {$\bullet$} at 3 4
\put {$\bullet$} at 3 3
\put {$\bullet$} at 3 2
\put {$\bullet$} at 3 1
\put {$\bullet$} at 3 0
\put {$\bullet$} at 3 -1
\put {$\star$} at 3 -2

\setlinear
\plot 4 3 3 2 /
\plot 3 4 3 3 /
\plot 3 3 3 2 /
\plot 3 2 3 1 /
\plot 3 1 3 0 /
\plot 3 0 3 -1 /
\plot 3 -1 3 -2 /

\put {$\bullet$} at 7 2
\put {$\bullet$} at 6 4
\put {$\bullet$} at 6 3
\put {$\bullet$} at 6 2
\put {$\bullet$} at 6 1
\put {$\bullet$} at 6 0
\put {$\bullet$} at 6 -1
\put {$\star$} at 6 -2

\setlinear
\plot 7 2 6 1 /
\plot 6 4 6 3 /
\plot 6 3 6 2 /
\plot 6 2 6 1 /
\plot 6 1 6 0 /
\plot 6 0 6 -1 /
\plot 6 -1 6 -2 /

\put {$\bullet$} at 10 1
\put {$\bullet$} at 9 4
\put {$\bullet$} at 9 3
\put {$\bullet$} at 9 2
\put {$\bullet$} at 9 1
\put {$\bullet$} at 9 0
\put {$\bullet$} at 9 -1
\put {$\star$} at 9 -2

\setlinear
\plot 10 1 9 0 /
\plot 9 4 9 3 /
\plot 9 3 9 2 /
\plot 9 2 9 1 /
\plot 9 1 9 0 /
\plot 9 0 9 -1 /
\plot 9 -1 9 -2 /

\put {$\bullet$} at 13 0
\put {$\bullet$} at 12 4
\put {$\bullet$} at 12 3
\put {$\bullet$} at 12 2
\put {$\bullet$} at 12 1
\put {$\bullet$} at 12 0
\put {$\bullet$} at 12 -1
\put {$\star$} at 12 -2

\setlinear
\plot 13 0 12 -1 /
\plot 12 4 12 3 /
\plot 12 3 12 2 /
\plot 12 2 12 1 /
\plot 12 1 12 0 /
\plot 12 0 12 -1 /
\plot 12 -1 12 -2 /

\endpicture \]

\[ \beginpicture
\setcoordinatesystem units <0.4cm,0.4cm>
\setplotarea x from 0 to 5, y from -1 to 1

\put {$\bullet$} at 0 3
\put {$\bullet$} at -1 2
\put {$\bullet$} at -2 3
\put {$\bullet$} at -2 2
\put {$\bullet$} at -2 1
\put {$\bullet$} at -2 0
\put {$\bullet$} at -2 -1
\put {$\star$} at -2 -2

\setlinear
\plot -2 1 -1 2 /
\plot -1 2 0 3 /
\plot -2 3 -2 2 /
\plot -2 2 -2 1 /
\plot -2 1 -2 0 /
\plot -2 0 -2 -1 /
\plot -2 -1 -2 -2 /

\put {$\bullet$} at 4 2
\put {$\bullet$} at 3 1
\put {$\bullet$} at 2 3
\put {$\bullet$} at 2 2
\put {$\bullet$} at 2 1
\put {$\bullet$} at 2 0
\put {$\bullet$} at 2 -1
\put {$\star$} at 2 -2

\setlinear
\plot 3 1 4 2 /
\plot 2 0 3 1 /
\plot 2 3 2 2 /
\plot 2 2 2 1 /
\plot 2 1 2 0 /
\plot 2 0 2 -1 /
\plot 2 -1 2 -2 /

\put {$\bullet$} at 8 1
\put {$\bullet$} at 7 0
\put {$\bullet$} at 6 3
\put {$\bullet$} at 6 2
\put {$\bullet$} at 6 1
\put {$\bullet$} at 6 0
\put {$\bullet$} at 6 -1
\put {$\star$} at 6 -2

\setlinear
\plot 7 0 8 1 /
\plot 6 -1 7 0 /
\plot 6 3 6 2 /
\plot 6 2 6 1 /
\plot 6 1 6 0 /
\plot 6 0 6 -1 /
\plot 6 -1 6 -2 /

\put {$\bullet$} at 13 2
\put {$\bullet$} at 12 1
\put {$\bullet$} at 11 0
\put {$\bullet$} at 10 2
\put {$\bullet$} at 10 1
\put {$\bullet$} at 10 0
\put {$\bullet$} at 10 -1
\put {$\star$} at 10 -2

\setlinear
\plot 11 0 12 1 /
\plot 12 1 13 2 /
\plot 10 -1 11 0 /
\plot 10 2 10 1 /
\plot 10 1 10 0 /
\plot 10 0 10 -1 /
\plot 10 -1 10 -2 /

\endpicture \]

\[ \beginpicture
\setcoordinatesystem units <0.4cm,0.4cm>
\setplotarea x from 0 to 5, y from -1 to 1

\put {$f=3$ [11]} at -7 0

\put {$\bullet$} at -2 3
\put {$\bullet$} at -2 2
\put {$\bullet$} at -3 3
\put {$\bullet$} at -3 2
\put {$\bullet$} at -3 1
\put {$\bullet$} at -3 0
\put {$\bullet$} at -3 -1
\put {$\star$} at -3 -2

\setlinear
\plot -2 2 -3 1 /
\plot -2 3 -3 2 /
\plot -3 3 -3 2 /
\plot -3 2 -3 1 /
\plot -3 1 -3 0 /
\plot -3 0 -3 -1 /
\plot -3 -1 -3 -2 /

\put {$\bullet$} at 1 3
\put {$\bullet$} at 1 1
\put {$\bullet$} at 0 3
\put {$\bullet$} at 0 2
\put {$\bullet$} at 0 1
\put {$\bullet$} at 0 0
\put {$\bullet$} at 0 -1
\put {$\star$} at 0 -2

\setlinear
\plot 0 0 1 1 /
\plot 0 2 1 3 /
\plot 0 3 0 2 /
\plot 0 2 0 1 /
\plot 0 1 0 0 /
\plot 0 0 0 -1 /
\plot 0 -1 0 -2 /

\put {$\bullet$} at 4 0
\put {$\bullet$} at 4 3
\put {$\bullet$} at 3 3
\put {$\bullet$} at 3 2
\put {$\bullet$} at 3 1
\put {$\bullet$} at 3 0
\put {$\bullet$} at 3 -1
\put {$\star$} at 3 -2

\setlinear
\plot 4 3 3 2 /
\plot 3 -1 4 0 /
\plot 3 3 3 2 /
\plot 3 2 3 1 /
\plot 3 1 3 0 /
\plot 3 0 3 -1 /
\plot 3 -1 3 -2 /

\put {$\bullet$} at 7 2
\put {$\bullet$} at 7 1
\put {$\bullet$} at 6 3
\put {$\bullet$} at 6 2
\put {$\bullet$} at 6 1
\put {$\bullet$} at 6 0
\put {$\bullet$} at 6 -1
\put {$\star$} at 6 -2

\setlinear
\plot 7 2 6 1 /
\plot 7 1 6 0 /
\plot 6 3 6 2 /
\plot 6 2 6 1 /
\plot 6 1 6 0 /
\plot 6 0 6 -1 /
\plot 6 -1 6 -2 /

\put {$\bullet$} at 10 0
\put {$\bullet$} at 10 2
\put {$\bullet$} at 9 3
\put {$\bullet$} at 9 2
\put {$\bullet$} at 9 1
\put {$\bullet$} at 9 0
\put {$\bullet$} at 9 -1
\put {$\star$} at 9 -2

\setlinear
\plot 10 0 9 -1 /
\plot 10 2 9 1 /
\plot 9 3 9 2 /
\plot 9 2 9 1 /
\plot 9 1 9 0 /
\plot 9 0 9 -1 /
\plot 9 -1 9 -2 /

\put {$\bullet$} at 13 0
\put {$\bullet$} at 13 1
\put {$\bullet$} at 12 3
\put {$\bullet$} at 12 2
\put {$\bullet$} at 12 1
\put {$\bullet$} at 12 0
\put {$\bullet$} at 12 -1
\put {$\star$} at 12 -2

\setlinear
\plot 13 0 12 -1 /
\plot 13 1 12 0 /
\plot 12 3 12 2 /
\plot 12 2 12 1 /
\plot 12 1 12 0 /
\plot 12 0 12 -1 /
\plot 12 -1 12 -2 /

\endpicture \]

\[ \beginpicture
\setcoordinatesystem units <0.4cm,0.4cm>
\setplotarea x from 0 to 5, y from -1 to 1

\put {$\bullet$} at 0 2
\put {$\bullet$} at -1 1
\put {$\bullet$} at -1 0
\put {$\bullet$} at -2 2
\put {$\bullet$} at -2 1
\put {$\bullet$} at -2 0
\put {$\bullet$} at -2 -1
\put {$\star$} at -2 -2

\setlinear
\plot 0 2 -1 1 /
\plot -1 1 -2 0 /
\plot -1 0 -2 -1 /
\plot -2 2 -2 1 /
\plot -2 1 -2 0 /
\plot -2 0 -2 -1 /
\plot -2 -1 -2 -2 /

\put {$\bullet$} at 3 0
\put {$\bullet$} at 4 1
\put {$\bullet$} at 3 1
\put {$\bullet$} at 2 2
\put {$\bullet$} at 2 1
\put {$\bullet$} at 2 0
\put {$\bullet$} at 2 -1
\put {$\star$} at 2 -2

\setlinear
\plot 3 0 2 -1 /
\plot 4 1 3 0 /
\plot 3 1 2 0 /
\plot 2 2 2 1 /
\plot 2 1 2 0 /
\plot 2 0 2 -1 /
\plot 2 -1 2 -2 /

\put {$\bullet$} at 8 1
\put {$\bullet$} at 7 0
\put {$\bullet$} at 7 2
\put {$\bullet$} at 6 2
\put {$\bullet$} at 6 1
\put {$\bullet$} at 6 0
\put {$\bullet$} at 6 -1
\put {$\star$} at 6 -2

\setlinear
\plot 7 0 8 1 /
\plot 6 -1 7 0 /
\plot 6 1 7 2 /
\plot 6 2 6 1 /
\plot 6 1 6 0 /
\plot 6 0 6 -1 /
\plot 6 -1 6 -2 /

\put {$\bullet$} at 12 2
\put {$\bullet$} at 11 2
\put {$\bullet$} at 11 1
\put {$\bullet$} at 10 2
\put {$\bullet$} at 10 1
\put {$\bullet$} at 10 0
\put {$\bullet$} at 10 -1
\put {$\star$} at 10 -2

\setlinear
\plot 11 1 12 2 /
\plot 11 1 11 2 /
\plot 10 0 11 1 /
\plot 10 2 10 1 /
\plot 10 1 10 0 /
\plot 10 0 10 -1 /
\plot 10 -1 10 -2 /

\put {$\bullet$} at 16 1
\put {$\bullet$} at 15 1
\put {$\bullet$} at 15 0
\put {$\bullet$} at 14 2
\put {$\bullet$} at 14 1
\put {$\bullet$} at 14 0
\put {$\bullet$} at 14 -1
\put {$\star$} at 14 -2

\setlinear
\plot 15 0 16 1 /
\plot 15 0 15 1 /
\plot 15 0 14 -1 /
\plot 14 2 14 1 /
\plot 14 1 14 0 /
\plot 14 0 14 -1 /
\plot 14 -1 14 -2 /

\endpicture \]

\[ \beginpicture
\setcoordinatesystem units <0.4cm,0.4cm>
\setplotarea x from 0 to 5, y from -1 to 1

\put {$f=4$ [2]} at -7 0

\put {$\bullet$} at 1.5 1
\put {$\bullet$} at 0.5 1
\put {$\bullet$} at -0.5 1
\put {$\bullet$} at -1.5 1
\put {$\bullet$} at 1 0
\put {$\bullet$} at -1 0
\put {$\bullet$} at 0 -1
\put {$\star$} at 0 -2

\setlinear
\plot 1.5 1 1 0 /
\plot 0.5 1 1 0 /
\plot -0.5 1 -1 0 /
\plot -1.5 1 -1 0 /
\plot 1 0 0 -1 /
\plot -1 0 0 -1 /
\plot 0 -1 0 -2 /

\put {$\bullet$} at 6 2
\put {$\bullet$} at 6 1
\put {$\bullet$} at 6 0
\put {$\bullet$} at 5 2
\put {$\bullet$} at 5 1
\put {$\bullet$} at 5 0
\put {$\bullet$} at 5 -1
\put {$\star$} at 5 -2

\setlinear
\plot 5 -2 5 -1 /
\plot 5 -1 5 0 /
\plot 5 0 5 1 /
\plot 5 1 5 2 /
\plot 5 -1 6 0 /
\plot 5 0 6 1 /
\plot 5 1 6 2 /

\endpicture \]

A concrete example is a permutation $G$, a 3-cycle in $\P_2^4$, such that
$$ G(1112)=1122, \;\;\; G(1122)=1222, \;\; \mbox{ and } \;\; G(1222)=1112. $$
The trees corresponding to $G$ are:
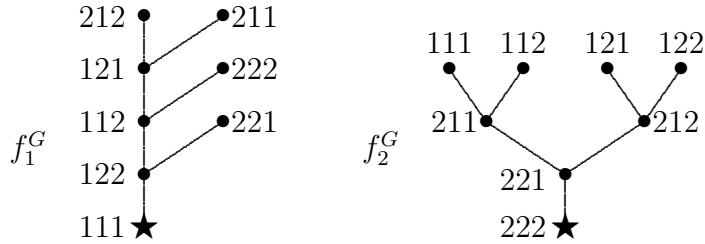
\[ \beginpicture
\setcoordinatesystem units <0.7cm,0.7cm>
\setplotarea x from 4 to 9, y from -1.5 to 3

\put {$\bullet$} at 5.8 2
\put {$\bullet$} at 7.2 2
\put {$\bullet$} at 4.2 2
\put {$\bullet$} at 2.8 2
\put {$\bullet$} at 3.5 1
\put {$\bullet$} at 6.5 1
\put {$\bullet$} at 5 0
\put {$\bigstar$} at 5 -1

\put {$122$} at 7.2 2.5
\put {$121$} at 5.8 2.5
\put {$111$} at 2.8 2.5
\put {$112$} at 4.2 2.5
\put {$211$} at 2.9 1
\put {$212$} at 7.1 1
\put {$221$} at 4.2 -0.1
\put {$222$} at 4.2 -1

\setlinear

\plot 7.2 2 6.5 1 /
\plot 5.8 2 6.5 1 /
\plot 2.8 2 3.5 1 /
\plot 4.2 2 3.5 1 /
\plot 3.5 1 5 0 /
\plot 6.5 1 5 0 /
\plot 5 0 5 -1 /

\put {$f_2^G$} at 1.5 0.5

\put {$\bullet$} at -1.5 3
\put {$\bullet$} at -1.5 2
\put {$\bullet$} at -1.5 1
\put {$\bullet$} at -3 3
\put {$\bullet$} at -3 2
\put {$\bullet$} at -3 1
\put {$\bullet$} at -3 0
\put {$\bigstar$} at -3 -1

\put {$212$} at -3.8 3
\put {$121$} at -3.8 2
\put {$112$} at -3.8 1
\put {$122$} at -3.8 0
\put {$111$} at -3.8 -1
\put {$211$} at -0.9 3
\put {$222$} at -0.9 2
\put {$221$} at -0.9 1

\setlinear

\plot -1.5 3 -3 2 /
\plot -1.5 2 -3 1 /
\plot -1.5 1 -3 0 /
\plot -3 3 -3 -1 /

\put {$f_1^G$} at -5.2 0.5

\endpicture \]
One checks that
\begin{equation}
\lambda_G^6={\rm id}
\end{equation}
but none of $\lambda_G$, $\lambda_G^2$, $\lambda_G^3$ is inner. Also note that
$\lambda_G(S_2)=S_2$.
}
\end{example}

\begin{example}
{\rm 
In the case of $\P_4^2$, we have 
\begin{align*}
& \# \{\sigma \in P_4^2 \ : \ \lambda_{u_\sigma}|_{\D_4} \in {\rm Aut}(\D_4) \} = 5400 \cdot 4!^4 = 1,791,590,400,\\
& \# \{\sigma \in P_4^2 \ : \ \lambda_{u_\sigma} \in {\rm Aut}(\O_4) \} = 5,771,520.
\end{align*}
In particular, there are $240,480$ distinct classes of automorphisms in ${\rm Out}({\mathcal O}_4)$
corresponding to permutations in $P_4^2$.

There are four (unlabeled) rooted trees with four vertices:
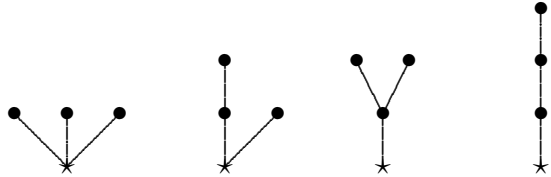
\[ \beginpicture
\setcoordinatesystem units <0.7cm,0.7cm>
\setplotarea x from 0 to 5, y from -1 to 0.1

\put {$\bullet$} at -4 0
\put {$\bullet$} at -3 0
\put {$\bullet$} at -2 0
\put {$\star$} at -3 -1

\setlinear
\plot -4 0 -3 -1 /
\plot -3 0 -3 -1 /
\plot -2 0 -3 -1 /

\put {$\bullet$} at 0 1
\put {$\bullet$} at 0 0
\put {$\bullet$} at 1 0
\put {$\star$} at 0 -1

\setlinear
\plot 0 1 0 0 /
\plot 1 0 0 -1 /
\plot 0 0 0 -1 /

\put {$\bullet$} at 3.5 1
\put {$\bullet$} at 2.5 1
\put {$\bullet$} at 3 0
\put {$\star$} at 3 -1

\setlinear
\plot 3.5 1 3 0 /
\plot 2.5 1 3 0 /
\plot 3 0 3 -1 /

\put {$\bullet$} at 6 2
\put {$\bullet$} at 6 1
\put {$\bullet$} at 6 0
\put {$\star$} at 6 -1

\setlinear
\plot 6 2 6 1 /
\plot 6 1 6 0 /
\plot 6 0 6 -1 /

\endpicture \]
and all of them arise from permutative automorphisms of $\O_4$. 
}
\end{example}

\bigskip

\noindent
Roberto Conti \\
Dipartimento di Scienze \\
Universit{\`a} di Chieti-Pescara `G. D'Annunzio' \\
Viale Pindaro 42, I--65127 Pescara, Italy \\
E-mail: conti@sci.unich.it \\

\smallskip \noindent
Wojciech Szyma{\'n}ski\\
Department of Mathematics and Computer Science \\
The University of Southern Denmark \\
Campusvej 55, DK-5230 Odense M, Denmark \\
E-mail: szymanski@imada.sdu.dk

\end{document}